\documentclass[10pt, a4paper]{article}
\usepackage{cmap}
\usepackage[main=russian,english]{babel}
\usepackage[utf8]{inputenc}
\usepackage[text={8in,10in},centering]{geometry}
\usepackage{amsfonts,amsmath,amsxtra,amsthm,amssymb,latexsym, tabularray}
\usepackage{graphicx}
\usepackage{verbatim}
\usepackage{cite, color}
\usepackage{subcaption}
\usepackage{float}
\usepackage{comment}
\theoremstyle{plain}
\newtheorem{theorem}{Теорема}[section]

\newtheorem{lemma}{\bf \indent Лемма}[section]
\newtheorem{corollary}{\bf \indent Следствие}[section]
\newtheorem*{theorem*}{\bf \indent Теорема}

\theoremstyle{remark}

\newtheorem{definition}{\bf \indent Определение}[section]
\newtheorem{example}{\bf \indent Пример}[section]

\numberwithin{equation}{section}

\newcommand{\be}[1]{\begin{equation}\label{#1}}
\newcommand{\ee}{\end{equation}}

\begin{document}
\parindent=0cm

\title{Секционная кривизна и структурные функции двумерного лоренцева многообразия}
\author{А.З. Али,\ Ю.Л. Сачков}
\date{\today}

\maketitle

\section{Введение}
В данной заметке мы получаем формулу для секционной кривизны на произвольном двумерном гладком многообразии $M$, снабжённом лоренцовой метрикой $g$ (симметрическое невырожденное гладкое тензорное поле второго порядка индекса один). Причём $g$ удовлетворяет следующему свойству: существуют базисные векторные поля $X_1$, $X_2$ на $M$ такие, что  $g(X_1,X_1) = a_{11}$, $g(X_2,X_2) = a_{22}$, $g(X_1,X_2) = a_{12} = g(X_2,X_1) = a_{21}$, где $a_{ij}$, $i,j=1,2$ --- постоянные вещественные числа.

В заметке приведена схема, по которой вычисляется секционная кривизна, вывод общей формулы, а также некоторые примеры применения полученного результата.

\section{Постановка задачи}
Напомним необходимые определения и известные утверждения \cite{beem}, \cite{lor_lob}. 

Пусть $M$ --- гладкое многообразие. Обозначим $Vec(M)$ множество гладких векторных полей на $M$. 

\begin{definition}
Связностью $D$ на многообразии $M$ называется отображение $D: (Vec(M))^2 \rightarrow Vec(M)$, удовлетворяющее аксиомам:
\begin{itemize}
    \item[$(1)$] $D_{fU+gV}W = fD_{U}W+gD_{V}W$ для любых $U$, $V$, $W$ $\in$ $Vec(M)$, для любых $f, g \in C^{\infty}(M)$,  
    \item[$(2)$] $D_{V}(\alpha_1 W_1 + \alpha_2 W_2) = \alpha_1 D_{V}W_1 + \alpha_2 D_{V}W_2$ для любых $W_1$, $W_2$ $\in Vec(M)$, для любых $\alpha_1$, $\alpha_2$ $\in \mathbb{R}$,
    \item[$(3)$] $D_{V}(fW) = (Vf)W + fD_{V}W$ для любых $V$, $W$ $\in Vec(M)$, для любой $f \in C^{\infty}(M)$.
\end{itemize}
\end{definition}

Векторное поле $D_{V}W$ называется {\it ковариантной производной поля $W$ вдоль $V$ для связности $D$}.

\begin{theorem}
На лоренцевом многообразии $(M, g)$ существует единственная связность $D$ такая, что:
\begin{itemize}
    \item[$(4)$] $[V, W] = D_{V}W - D_{W}V$,
    \item[$(5)$] $Xg(V, W) = g(D_X V, W) + g(V, D_{X} W)$, 
\end{itemize}
для любых $X, V, W \in Vec(M)$. $D$ называется {\it связностью Леви-Чивиты на $M$}, и характеризуется формулой Кошуля
\begin{equation}
\label{koszul}
2g(D_{V}W,X) =Vg(W,X)+Wg(X,V)-Xg(V,W)-g(V,[W,X])+g(W,[X,V])+g(X,[V,W]).
\end{equation}
\end{theorem}

Отображение $R: (Vec(M))^3 \rightarrow C^{\infty}(M)$, задаваемое формулой $R_{XY}Z = D_{[X,Y]}Z-[D_{X},D_{Y}]Z$, где $X, Y, Z \in Vec(M)$, называется {\it тензором кривизны Римана $(M,g)$}.

Пусть $q \in M$, $P$ - двумерная плоскость в касательном пространстве $T_{q}M$. Тогда для векторов $v, w \in T_{q}M$ определим $Q(v,w) = g(v,v)g(w,w)-(g(v,w))^2$. Плоскость $P$ называется {\it невырожденной}, если $Q(v,w) \neq 0$ для некоторого (а тогда и любого) базиса $v, w$ плоскости $P$.

\begin{lemma}
Пусть $P \subseteq T_{q}M$ - невырожденная двумерная плоскость. Тогда число
\begin{equation}
\label{sectionalcurvature}
    K(q,P) = \frac{g(R_{vw}v,w)}{Q(v,w)}
\end{equation}
не зависит от выбора базиса в плоскости $P$, и называется {\bf секционной кривизной плоского сектора $P$}.
\end{lemma}

\begin{theorem}
Пусть $(M,g)$ - лоренцево многообразие размерности $n \geqslant 2$, и пусть $K \in \mathbb{R}$. Тогда следующие условия эквивалентны:
\begin{itemize}
    \item[$(1)$] $M$ имеет постоянную кривизну $K$, 
    \item[$(2)$] для любой точки $q \in M$ существует окрестность, изометричная открытому подмножеству пространства де-Ситтера $\mathbb{S}_1^n$ для $K>0$, пространству Минковского $\mathbb{R}_{1}^{n}$ для $K = 0$ или пространству анти-де Ситтера $\widetilde{\mathbb{H}_{1}^{n}}$ для $K < 0$. 
\end{itemize}
\end{theorem}

\section{План вычисления секционной кривизны для двумерного многообразия}

\begin{itemize}
    \item[$(1)$] Вычислить коммутатор $[X_1, X_2]$;
    \item[$(2)$] Вычислить коэффициенты связности Леви-Чивиты с помощью формулы Кошуля (\ref{koszul}), используя базисные векторные поля $X_1, X_2$;
    \item[$(3)$] Вычислить тензор кривизны Римана $R_{v w}v$, где $v = X_1$, $w = X_2$;
    \item[$(4)$] Вычислить $Q(v,w)$, где $v = X_1$, $w = X_2$;
    \item[$(5)$] Получить секционную кривизну $K$, используя проделанные вычисления и формулу (\ref{sectionalcurvature}).
\end{itemize}

\section{Основная формула, полезные следствия и примеры}

\begin{theorem}
Секционная кривизна лоренцева двумерного многообразия $M$ с лоренцевой метрикой $g$, её базисными векторными полями $X_1$, $X_2$, $[X_1,X_2] = c_{12}^{1}X_1+c_{12}^{2}X_2$ и постоянными значениями метрики на них: $g(X_1,X_1) = a_{11}$, $g(X_2,X_2) = a_{22}$, $g(X_1,X_2) = a_{12} = g(X_2,X_1) = a_{21}$, равна 
\begin{equation}
\label{dvum_kriv_obsh}
    K= \frac{a_{12}\xi_1 + a_{22}\xi_2}{\det{\mathcal{A}}^2},
\end{equation}
где 
\begin{equation*}
    \xi_1 = a_{12} \left[ \left( c_{12}^1A + c_{12}^2B \right) - \left( X_1 B \right) + \left( X_2 A \right) \right], \quad \xi_2 = a_{11} \left[ - \left( c_{12}^1A + c_{12}^2B \right) + \left( X_1B \right) - \left( X_2A \right) \right],
\end{equation*}
\begin{equation*}
A = a_{11}c_{12}^1+a_{12}c_{12}^2, \quad B = a_{12}c_{12}^1+a_{22}c_{12}^2,
\end{equation*}

\begin{equation}
\label{dvum_matr}
    \mathcal{A} = 
    \begin{pmatrix}
    a_{11} & a_{12}\\
    a_{12} & a_{22}
    \end{pmatrix}. 
\end{equation}

\end{theorem}

\begin{corollary}
Для ортогональных $X_1$, $X_2$, т.е. $a_{12} = a_{21} = 0$, $a_{11} \neq 0$, $a_{22} \neq 0$, секционная кривизна равна
\begin{equation}
\label{dvum_kriv_ortogon}
K = \frac{-a_{11}\left(c_{12}^1\right)^2 -a_{22} \left(c_{12}^2\right)^2 - a_{11} \left(X_2 c_{12}^1\right) + a_{22} \left(X_1 c_{12}^2\right)}{a_{11}}.
\end{equation}
\end{corollary}

\begin{corollary}
Для ортонормированных $X_1$, $X_2$, т.е. $a_{11} = -1$, $a_{12} = a_{21} = 0$, $a_{22} = 1$, секционная кривизна равна
\begin{equation}
\label{dvum_kriv_ortonorm}
K = -\left(c_{12}^1\right)^2 + \left(c_{12}^2\right)^2 - \left(X_2 c_{12}^1\right) - \left(X_1 c_{12}^2\right).
\end{equation}
\end{corollary}

\begin{example}
{\bf Двумерное пространство анти-де Ситтера $\widetilde{\mathbb{H}_{1}^{2}}$.}

Лоренцева метрика имеет выражение: $g =  - \ch^2{\theta} d \varphi^2+d \theta^2$, и ортонормированный репер: $X_1 = \frac{1}{\ch{\theta}} \frac{\partial}{\partial \varphi},\ X_2 = \frac{\partial }{\partial \theta}$. Имеем $g(X_1,X_1) = -1$, $g(X_2,X_2) = 1$, $g(X_1,X_2) = g(X_2,X_1) = 0$, а также $[X_1,X_2] = \th{\theta}X_1$, то есть структурные функции следующие: $c_{12}^1 = \th{\theta}$, $c_{12}^2 = 0$. 

Поэтому, по формуле (\ref{dvum_kriv_ortonorm}), секционная кривизна двумерного пространства анти-де Ситтера равна:
\begin{equation*}
K = -\left[\left(c_{12}^1\right)^2 - \left(c_{12}^2\right)^2 + \left(X_2 c_{12}^1\right) + \left(X_1 c_{12}^2\right)\right] = -\left[ (\th{\theta})^2 + \partial_{\theta}(\th{\theta}) \right] = -\frac{\sh^2{\theta}+1}{\ch^2{\theta}} = -1,
\end{equation*}
что совпадает с полученным в \cite{lor_lob}.
\end{example}

\begin{example}
\label{igreki}
Пусть $X_1$, $X_2$ --- базис, т.ч. $[X_1,X_2] = -X_1$, а $Y_1$, $Y_2$ --- ортонормированный базис лоренцевой метрики $g$, причём $Y_1 = \alpha X_1 + \gamma X_2$, $Y_2 = \beta X_1 + \delta X_2$. 

Посчитаем секционную кривизну по формуле (\ref{dvum_kriv_ortonorm}). Для этого требуется сперва посчитать структурные функции.
\begin{equation*}
[Y_1,Y_2] = [\alpha X_1 + \gamma X_2,\beta X_1 + \delta X_2] = \alpha \beta[X_1,X_1] + \alpha \delta [X_1,X_2] + \gamma \beta [X_2,X_1] + \gamma \delta [X_2,X_2] = (\gamma \beta - \alpha \delta)X_1.
\end{equation*}
Разложим полученный вектор по $Y_1$, $Y_2$:
\begin{equation*}
\gamma Y_2 - \delta Y_1 = \gamma( \beta X_1 + \delta X_2 ) - \delta ( \alpha X_1 + \gamma X_2 ) = (\gamma \beta - \alpha \delta)X_1,
\end{equation*}
поэтому
\begin{equation*}
[Y_1,Y_2] = \gamma Y_2 - \delta Y_1 = c_{12}^1 Y_1 + c_{12}^2 Y_2,
\end{equation*}
то есть $c_{12}^1 = -\delta$, $c_{12}^2 = \gamma$.
Получаем:
\begin{equation*}
K = - \left[ \left(c_{12}^1\right)^2 - \left(c_{12}^2\right)^2 + \left(Y_2 c_{12}^1\right) + \left(Y_1 c_{12}^2\right) \right] = - \left[ \delta^2 - \gamma^2 \right] = \gamma^2 - \delta^2,
\end{equation*}
что совпадает с полученным в \cite{lor_lob}.
\end{example}

\begin{example}
{\bf Левоинвариантные задачи на двумерной группе Ли $G$.}

Пусть $X_1$, $X_2$ --- векторные поля, такие, что $[X_2,X_1] = X_1$, то есть $c_{12}^1 = -1$, $c_{12}^2 = 0$ в наших обозначениях, а $g(X_1,X_1) = c^2-a^2 = a_{11}$, $g(X_1,X_2) = g(X_2, X_1) = cd-ab = a_{12}$, $g(X_2,X_2) = d^2-b^2 = a_{22}$. 

Вычислим секционную кривизну по формуле (\ref{dvum_kriv_ortogon}). 

Так как структурные функции являются константами, все производные Ли вдоль наших векторных полей равны нулю. Поэтому достаточно посчитать $A$ и $B$, чтобы вычислить $\xi_1$, $\xi_2$.
\begin{equation*}
A = a_{11}c_{12}^1 + a_{12}c_{12}^2 = -a_{11} = a^2-c^2,\ B = a_{12}c_{12}^1+a_{22}c_{12}^2 = -a_{12}= ab-cd,
\end{equation*}
поэтому 
\begin{equation*}
\xi_1 = a_{12} \left[ \left( c_{12}^1A + c_{12}^2B \right) - \left( X_1 B \right) + \left( X_2 A \right) \right] = -a_{12}A = -a_{12}( a_{11}c_{12}^1 + a_{12}c_{12}^2 ) = a_{11}a_{12},
\end{equation*}
\begin{equation*}
\xi_2 =  a_{11} \left[ - \left( c_{12}^1A + c_{12}^2B \right) + \left( X_1B \right) - \left( X_2A \right) \right] = a_{11}( a_{11}c_{12}^1 + a_{12}c_{12}^2) = -(a_{11})^2.
\end{equation*}
\begin{equation*}
a_{12}\xi_1 + a_{22}\xi_2 = a_{12}a_{11}a_{12} - a_{22}(a_{11})^2 = -a_{11}\left( a_{11}a_{22} - (a_{12})^2 \right) = -a_{11}\det{\mathcal{A}}.
\end{equation*}
Также отметим, что
\begin{equation*}
\det{\mathcal{A}} = a_{11}a_{22}-(a_{12})^2 = (c^2-a^2)(d^2-b^2)-(cd-ab)^2 = -(ad-bc)^2 = -\det(\mathcal{B})^2, 
\end{equation*}
где 
\begin{equation*}
\mathcal{B} = 
\begin{pmatrix}
a & b \\
c & d
\end{pmatrix}.
\end{equation*}
Тогда
\begin{equation*}
K = \frac{a_{12}\xi_1 + a_{22}\xi_2}{\det{\mathcal{A}}^2} = \frac{-a_{11}\det{\mathcal{A}}}{\det(\mathcal{B})^4} = \frac{a_{11}}{\det(\mathcal{B})^2},
\end{equation*}

что совпадает с полученным в \cite{lor_lob}.

\end{example}

\section{Вывод основной формулы}

Действуем по плану.

\begin{itemize}
    \item[$(1)$] Коммутатор мы записали в общем виде:
    \begin{equation*}
        [X_1,X_2] = c_{12}^{1}X_1+c_{12}^{2}X_2.
    \end{equation*}
    \item[$(2)$] Так как для $X_1, X_2$ $a_{ij} = const$, $\forall$ $i,j = 1,2$, первые три члена в правой части формулы Кошуля (\ref{koszul}) равны нулю, если вместо $V, W, X$ подставлять $X_1, X_2$ в любых комбинациях. Поэтому имеют смысл только четвёртый, пятый и шестой элементы правой части формулы (\ref{koszul}). Всего требуется вычислить 8 коэффициентов, а именно $g(D_{X_1}X_1,X_i),\ i=1,2$; $g(D_{X_1}X_2,X_i),\ i=1,2$; $g(D_{X_2}X_1,X_i),\ i=1,2$; $g(D_{X_2}X_2,X_i),\ i=1,2$. Приступим к вычислениям. Сразу можно сказать, что $g(D_{X_1}X_1,X_1) = 0$, $g(D_{X_2}X_2,X_2) = 0$, так как все коммутаторы состоят из, соответственно, $[X_1,X_1]$ или $[X_2,X_2]$, равных нулю тождественно. 
    \begin{equation*}
    g(D_{X_1}X_1,X_2) = \frac{1}{2}(-g(X_1,[X_1,X_2])+g(X_1,[X_2,X_1])+g(X_2,[X_1,X_1]))= 
    \end{equation*}
    \begin{equation*}
    = \frac{1}{2}( -g(X_1,c_{12}^{1}X_1+c_{12}^{2}X_2) + g(X_1,-c_{12}^{1}X_1-c_{12}^{2}X_2)) =
    \end{equation*}
    \begin{equation*}
    = \frac{1}{2}( -a_{11}c_{12}^1-a_{12}c_{12}^2-a_{11}c_{12}^1-a_{12}c_{12}^2 ) = -a_{11}c_{12}^1-a_{12}c_{12}^2 = -A,
    \end{equation*}
    \begin{equation*}
    g(D_{X_1}X_2,X_1) = \frac{1}{2}(-g(X_1,[X_2,X_1])+g(X_2,[X_1,X_1])+g(X_1,[X_1,X_2])) = 
    \end{equation*}
    \begin{equation*}
    = \frac{1}{2}(-g(X_1,-c_{12}^{1}X_1-c_{12}^{2}X_2)+g(X_1,c_{12}^{1}X_1+c_{12}^{2}X_2)) =
    \end{equation*}
    \begin{equation*}
    = \frac{1}{2}( a_{11}c_{12}^1 + a_{12}c_{12}^2 + a_{11}c_{12}^1 + a_{12}c_{12}^2 ) = a_{11}c_{12}^1 + a_{12}c_{12}^2 = A,
    \end{equation*}
    \begin{equation*}
    g(D_{X_1}X_2,X_2) = \frac{1}{2}( -g(X_1,[X_2,X_2])+g(X_2,[X_2,X_1])+g(X_2,[X_1,X_2]) ) = \end{equation*}
    \begin{equation*}
    = \frac{1}{2}( g(X_2,-c_{12}^{1}X_1-c_{12}^{2}X_2)+g(X_2,c_{12}^{1}X_1+c_{12}^{2}X_2) ) = 0,
    \end{equation*}
    \begin{equation*}
    g(D_{X_2}X_1,X_1) = \frac{1}{2}( -g(X_2,[X_1,X_1])+g(X_1,[X_1,X_2])+g(X_1,[X_2,X_1]) ) =
    \end{equation*}
    \begin{equation*}
    =  \frac{1}{2}( g(X_1,c_{12}^{1}X_1+c_{12}^{2}X_2) + g(X_1,-c_{12}^{1}X_1-c_{12}^{2}X_2) ) = 0,
    \end{equation*}
    \begin{equation*} g(D_{X_2}X_1,X_2) = \frac{1}{2}( -g(X_2,[X_1,X_2])+g(X_1,[X_2,X_2])+g(X_2,[X_2,X_1]) ) =
    \end{equation*}
    \begin{equation*}
    = \frac{1}{2}( -g(X_2,c_{12}^{1}X_1+c_{12}^{2}X_2) + g(X_2,-c_{12}^{1}X_1-c_{12}^{2}X_2) ) =
    \end{equation*}
    \begin{equation*} 
    = \frac{1}{2}( -a_{12}c_{12}^1-a_{22}c_{12}^2-a_{12}c_{12}^1 - a_{22}c_{12}^2 ) = -a_{12}c_{12}^1-a_{22}c_{12}^2 = -B,
    \end{equation*}
    \begin{equation*} g(D_{X_2}X_2,X_1) = \frac{1}{2}( -g(X_2,[X_2,X_1])+g(X_2,[X_1,X_2])+g(X_1,[X_2,X_2]) ) = 
    \end{equation*}
    \begin{equation*}
    = \frac{1}{2}( -g(X_2,-c_{12}^{1}X_1-c_{12}^{2}X_2)+g(X_2,c_{12}^{1}X_1+c_{12}^{2}X_2) ) 
    \end{equation*}
    \begin{equation*} 
    = \frac{1}{2}( a_{12}c_{12}^1+a_{22}c_{12}^2 + a_{12}c_{12}^1 + a_{22}c_{12}^2 ) = a_{12}c_{12}^1+a_{22}c_{12}^2 = B.
    \end{equation*}

    Теперь вычислим ковариантные производные.

    \begin{itemize}
        \item Временно обозначим 
        \begin{equation*}
            D_{X_1}X_1 = uX_1 + vX_2.
        \end{equation*}
        Найдём $u$ и $v$ из системы уравнений:
        \begin{equation*}
        \begin{cases}
            g\left( D_{X_1}X_1,X_1 \right) = 0,\\
            g\left(D_{X_1}X_1,X_2 \right) = -A.
        \end{cases}
        \end{equation*}
        \begin{equation*}
        \begin{cases}
            g\left( uX_1 + vX_2,X_1 \right) = 0,\\
            g\left(uX_1 + vX_2,X_2 \right) = -A,
        \end{cases}
        \Leftrightarrow
        \begin{cases}
    ug\left(X_1,X_1\right)+vg\left(X_2,X_1 \right) = 0,\\
    ug\left(X_1,X_2\right)+vg\left(X_2,X_2 \right) = -A,
        \end{cases}
        \Leftrightarrow
        \begin{cases}
            ua_{11} + va_{12} = 0,\\
            ua_{12} + va_{22} = -A,
        \end{cases}
        \Leftrightarrow
        \end{equation*}
    \begin{equation*}
    \Leftrightarrow
    \begin{cases}
        u = \frac{a_{12}A}{a_{11}a_{22}-(a_{12})^2} = \frac{a_{12}A}{\det{(\mathcal{A})}},\\
        v = \frac{-Aa_{11}}{a_{11}a_{22}-(a_{12})^2} = \frac{-Aa_{11}}{\det{(\mathcal{A})}}.
    \end{cases}
    \end{equation*}

    \item Временно обозначим 
        \begin{equation*}
            D_{X_1}X_2 = uX_1 + vX_2.
        \end{equation*}
        Найдём $u$ и $v$ из системы уравнений:
        \begin{equation*}
        \begin{cases}
            g\left( D_{X_1}X_2,X_1 \right) = A,\\
            g\left(D_{X_1}X_2,X_2 \right) = 0.
        \end{cases}
        \end{equation*}
        \begin{equation*}
        \begin{cases}
        g\left( uX_1 + vX_2,X_1 \right) = A,\\
        g\left(uX_1 + vX_2,X_2 \right) = 0.
        \end{cases}
        \Leftrightarrow
        \begin{cases}
        ua_{11} + va_{12} = A,\\
        ua_{12} + va_{22} = 0,
        \end{cases}
        \Leftrightarrow
        \begin{cases}
        u = \frac{a_{22}A}{\det{(\mathcal{A})}},\\
        v = \frac{-a_{12}A}{\det{(\mathcal{A})}}.
        \end{cases}
        \end{equation*}

        \item Временно обозначим 
        \begin{equation*}
            D_{X_2}X_1 = uX_1 + vX_2.
        \end{equation*}
        Найдём $u$ и $v$ из системы уравнений:
        \begin{equation*}
        \begin{cases}
            g\left( D_{X_2}X_1,X_1 \right) = 0,\\
            g\left(D_{X_2}X_1,X_2 \right) = -B.
        \end{cases}
        \end{equation*}
        \begin{equation*}
        \begin{cases}
        g\left( uX_1 + vX_2,X_1 \right) = 0,\\
            g\left(uX_1 + vX_2,X_2 \right) = -B,
        \end{cases}
        \Leftrightarrow
        \begin{cases}
        ua_{11} + va_{12} = 0,\\
        ua_{12} + va_{22} = -B,
        \end{cases}
        \Leftrightarrow
        \begin{cases}
        u = \frac{a_{12}B}{\det{(\mathcal{A})}},\\
        v = \frac{-a_{11}B}{\det{(\mathcal{A})}}.
        \end{cases}
        \end{equation*}

        \item Временно обозначим 
        \begin{equation*}
            D_{X_2}X_2 = uX_1 + vX_2.
        \end{equation*}
        Найдём $u$ и $v$ из системы уравнений:
        \begin{equation*}
        \begin{cases}
            g\left( D_{X_2}X_2,X_1 \right) = B,\\
            g\left(D_{X_2}X_2,X_2 \right) = 0.
        \end{cases}
        \end{equation*}
        \begin{equation*}
        \begin{cases}
            g\left( uX_1 + vX_2,X_1 \right) = B,\\
            g\left(uX_1 + vX_2,X_2 \right) = 0.
        \end{cases}
        \Leftrightarrow
            \begin{cases}
                ua_{11} + va_{12} = B,\\
                ua_{12} + va_{22} = 0,
            \end{cases}
        \Leftrightarrow
        \begin{cases}
            u = \frac{a_{22}B}{\det{(\mathcal{A})}},\\
            v = \frac{-a_{12}B}{\det{(\mathcal{A})}}.
        \end{cases}
        \end{equation*}
    \end{itemize}
    
    Поэтому
    \begin{align*}
    & D_{X_1}X_1 = \frac{A}{\det{(\mathcal{A})}}\left( a_{12}X_1 - a_{11}X_2 \right),\\
    & D_{X_1}X_2 = \frac{A}{\det{(\mathcal{A})}}\left( a_{22}X_1 - a_{12}X_2 \right),\\
    & D_{X_2}X_1 = \frac{B}{\det{(\mathcal{A})}}\left( a_{12}X_1 - a_{11}X_2 \right),\\
    & D_{X_2}X_2 = \frac{B}{\det{(\mathcal{A})}}\left( a_{22}X_1 - a_{12}X_2 \right).
    \end{align*}
    \item[$(3)$] \begin{align*}
    & R_{X_1 X_2}X_1 = D_{[X_1,X_2]}X_1 - (D_{X_1}D_{X_2}-D_{X_2}D_{X_1})X_1 = \\  
    & = D_{c_{12}^{1}X_1+c_{12}^{2}X_2}X_1 - D_{X_1}\frac{B}{\det{(\mathcal{A})}}\left( a_{12}X_1 - a_{11}X_2 \right) + D_{X_2}\frac{A}{\det{(\mathcal{A})}}\left( a_{12}X_1 - a_{11}X_2 \right) = \\
    & = D_{c_{12}^{1}X_1+c_{12}^{2}X_2}X_1 - \frac{a_{12}}{\det{(\mathcal{A})}}D_{X_1}\left( BX_1 \right) + \frac{a_{11}}{\det{(\mathcal{A})}}D_{X_1}\left( BX_2 \right) + \frac{a_{12}}{\det{(\mathcal{A})}}D_{X_2}\left( A X_1 \right) - \frac{a_{11}}{\det{(\mathcal{A})}}D_{X_2}\left( A X_2 \right).
    \end{align*}
    
    Посчитаем отдельно каждое слагаемое, а затем соберём в одну формулу.

    \begin{align*}
    & D_{c_{12}^{1}X_1+c_{12}^{2}X_2}X_1 = c_{12}^1D_{X_1}X_1 + c_{12}^2D_{X_2}X_1 = c_{12}^1\frac{A}{\det{(\mathcal{A})}}\left( a_{12}X_1 - a_{11}X_2 \right) + c_{12}^2\frac{B}{\det{(\mathcal{A})}}\left( a_{12}X_1 - a_{11}X_2 \right) = \\
    & = \frac{1}{\det{(\mathcal{A})}}\left( a_{12}\left( c_{12}^1A + c_{12}^2B \right)X_1 - a_{11}\left( c_{12}^1A + c_{12}^2B \right)X_2 \right).
    \end{align*}

    \begin{align*}
    & D_{X_1}\left( B X_1 \right) = D_{X_1}\left( \left( a_{12}c_{12}^1+a_{22}c_{12}^2 \right) X_1 \right) = a_{12}\left[ (X_1c_{12}^1)X_1 + c_{12}^1D_{X_1}X_1 \right] + a_{22}\left[ (X_1 c_{12}^2)X_1 + c_{12}^2D_{X_1}X_1\right] = \\
    & =  a_{12}(X_1c_{12}^1)X_1 + a_{22}(X_1 c_{12}^2)X_1 + a_{12}c_{12}^1\frac{A}{\det{(\mathcal{A})}}\left( a_{12}X_1 - a_{11}X_2 \right) + a_{22}c_{12}^2\frac{A}{\det{(\mathcal{A})}}\left( a_{12}X_1 - a_{11}X_2 \right)=\\
    & = \left[ a_{12}(X_1c_{12}^1) + a_{22}(X_1 c_{12}^2) + \frac{A}{\det{(\mathcal{A})}}(c_{12}^1(a_{12})^2 + c_{12}^2a_{12}a_{22}) \right]X_1 - \frac{A}{\det{(\mathcal{A})}}\left[ c_{12}^1a_{11}a_{12} + c_{12}^2a_{11}a_{22} \right]X_2 = \\
    & = \left[ a_{12}(X_1c_{12}^1) + a_{22}(X_1 c_{12}^2) + \frac{a_{12}AB}{\det{(\mathcal{A})}} \right]X_1 - \frac{a_{11}AB}{\det{(\mathcal{A})}}X_2.
    \end{align*}

    \begin{align*}
    & D_{X_1}\left( B X_2 \right) = D_{X_1}\left( \left( a_{12}c_{12}^1+a_{22}c_{12}^2 \right) X_2 \right) = a_{12}\left[ \left( X_1c_{12}^1 \right)X_2 + c_{12}^1D_{X_1}X_2 \right] + a_{22}\left[ \left( X_1c_{12}^2 \right)X_2 + c_{12}^2D_{X_1}X_2 \right] = \\
    & = a_{12}\left( X_1c_{12}^1 \right)X_2 + a_{22}\left( X_1c_{12}^2 \right)X_2 + a_{12}c_{12}^1\frac{A}{\det{(\mathcal{A})}}\left( a_{22}X_1 - a_{12}X_2 \right) + a_{22}c_{12}^2\frac{A}{\det{(\mathcal{A})}}\left( a_{22}X_1 - a_{12}X_2 \right)=\\
    & = \frac{A}{\det{(\mathcal{A})}}\left[ c_{12}^1a_{12}a_{22} + c_{12}^2\left(a_{22}\right)^2 \right]X_1 + \left[ a_{12}\left( X_1c_{12}^1 \right) + a_{22}\left( X_1c_{12}^2 \right) - \frac{A}{\det{(\mathcal{A})}}\left( c_{12}^1\left(a_{12}\right)^2 + c_{12}^2a_{12}a_{22} \right) \right]X_2 = \\
    & = \frac{a_{22}AB}{\det{(\mathcal{A})}}X_1 + \left[ a_{12}\left( X_1c_{12}^1 \right) + a_{22}\left( X_1c_{12}^2 \right) - \frac{a_{12}AB}{\det{(\mathcal{A})}} \right]X_2.
    \end{align*}

    \begin{align*}
    & D_{X_2}\left( A X_1 \right) = D_{X_2}\left( \left( a_{11}c_{12}^1 + a_{12}c_{12}^2 \right) X_1 \right) = a_{11}\left[ \left( X_2 c_{12}^1 \right)X_1 + c_{12}^1D_{X_2}X_1 \right] + a_{12}\left[ \left( X_2 c_{12}^2 \right)X_1 + c_{12}^2D_{X_2}X_1 \right] = \\
    & = a_{11}\left( X_2 c_{12}^1 \right)X_1 + a_{12}\left( X_2 c_{12}^2 \right)X_1 + a_{11}c_{12}^1\frac{B}{\det{(\mathcal{A})}}\left( a_{12}X_1 - a_{11}X_2 \right) + a_{12}c_{12}^2\frac{B}{\det{(\mathcal{A})}}\left( a_{12}X_1 - a_{11}X_2 \right) = \\
    & = \left[ a_{11}\left( X_2 c_{12}^1 \right) + a_{12}\left( X_2 c_{12}^2 \right) + \frac{B}{\det{(\mathcal{A})}}\left( c_{12}^1a_{11}a_{12} + c_{12}^2\left( a_{12}\right)^2 \right)  \right]X_1 - \frac{B}{\det{(\mathcal{A})}}\left[ c_{12}^1\left( a_{11} \right)^2 + c_{12}^2a_{11}a_{12} \right]X_2 = \\
    & = \left[ a_{11}\left( X_2 c_{12}^1 \right) + a_{12}\left( X_2 c_{12}^2 \right) + \frac{a_{12}AB}{\det{(\mathcal{A})}}  \right]X_1 - \frac{a_{11}AB}{\det{(\mathcal{A})}}X_2.
    \end{align*}
    
    \begin{align*}
    & D_{X_2}\left( A X_2 \right) = D_{X_2}\left( \left( a_{11}c_{12}^1 + a_{12}c_{12}^2 \right)X_2 \right) = a_{11}\left[ \left( X_2c_{12}^1 \right)X_2 + c_{12}^1D_{X_2}X_2 \right] + a_{12}\left[ \left( X_2c_{12}^2 \right)X_2 
 +c_{12}^2D_{X_2}X_2 \right] = \\
 & = a_{11}\left( X_2c_{12}^1 \right)X_2 + a_{12}\left( X_2c_{12}^2 \right)X_2 + a_{11}c_{12}^1\frac{B}{\det{(\mathcal{A})}}\left( a_{22}X_1 - a_{12}X_2 \right) + a_{12}c_{12}^2\frac{B}{\det{(\mathcal{A})}}\left( a_{22}X_1 - a_{12}X_2 \right) = \\
 & = \frac{B}{\det{(\mathcal{A})}}\left[ c_{12}^1a_{11}a_{22} + c_{12}^2a_{12}a_{22} \right]X_1 + \left[ a_{11}\left( X_2c_{12}^1 \right) + a_{12}\left( X_2c_{12}^2 \right) - \frac{B}{\det{(\mathcal{A})}}\left( c_{12}^1a_{11}a_{12} + c_{12}^2 \left(a_{12}\right)^2 \right) \right]X_2 = \\
 & = \frac{a_{22}AB}{\det{(\mathcal{A})}} X_1 + \left[ a_{11}\left( X_2c_{12}^1 \right) + a_{12}\left( X_2c_{12}^2 \right) - \frac{a_{12}AB}{\det{(\mathcal{A})}} \right]X_2.
    \end{align*}

    Получаем

    \begin{align*}
    &  D_{c_{12}^{1}X_1+c_{12}^{2}X_2}X_1 - \frac{a_{12}}{\det{(\mathcal{A})}}D_{X_1}\left( BX_1 \right) + \frac{a_{11}}{\det{(\mathcal{A})}}D_{X_1}\left( BX_2 \right) + \frac{a_{12}}{\det{(\mathcal{A})}}D_{X_2}\left( A X_1 \right) - \frac{a_{11}}{\det{(\mathcal{A})}}D_{X_2}\left( A X_2 \right) = \\
    & = \frac{1}{\det{(\mathcal{A})}}\left( a_{12}\left( c_{12}^1A + c_{12}^2B \right)X_1 - a_{11}\left( c_{12}^1A + c_{12}^2B \right)X_2 \right) - \\
    & - \frac{a_{12}}{\det{(\mathcal{A})}}\left( \left[ a_{12}(X_1c_{12}^1) + a_{22}(X_1 c_{12}^2) + \frac{a_{12}AB}{\det{(\mathcal{A})}} \right]X_1 - \frac{a_{11}AB}{\det{(\mathcal{A})}}X_2 \right)+\\
    & + \frac{a_{11}}{\det{(\mathcal{A})}}\left( \frac{a_{22}AB}{\det{(\mathcal{A})}}X_1 + \left[ a_{12}\left( X_1c_{12}^1 \right) + a_{22}\left( X_1c_{12}^2 \right) - \frac{a_{12}AB}{\det{(\mathcal{A})}} \right]X_2 \right) +\\
    & + \frac{a_{12}}{\det{(\mathcal{A})}}\left( \left[ a_{11}\left( X_2 c_{12}^1 \right) + a_{12}\left( X_2 c_{12}^2 \right) + \frac{a_{12}AB}{\det{(\mathcal{A})}}  \right]X_1 - \frac{a_{11}AB}{\det{(\mathcal{A})}}X_2 \right) - \\
    & - \frac{a_{11}}{\det{(\mathcal{A})}}\left( \frac{a_{22}AB}{\det{(\mathcal{A})}} X_1 + \left[ a_{11}\left( X_2c_{12}^1 \right) + a_{12}\left( X_2c_{12}^2 \right) - \frac{a_{12}AB}{\det{(\mathcal{A})}} \right]X_2 \right) = \\
    & = \frac{1}{\det{\mathcal{A}}}\left(a_{12} \left[ \left( c_{12}^1A + c_{12}^2B \right) - \left( X_1 B \right) + \left( X_2 A \right) \right]X_1 + a_{11} \left[ - \left( c_{12}^1A + c_{12}^2B \right) + \left( X_1B \right) - \left( X_2A \right) \right]X_2 \right) = \frac{1}{\det{\mathcal{A}}}\left( \xi_1 X_1 + \xi_2 X_2 \right).
    \end{align*}
    \item[$(4)$]
    \begin{equation*}
    Q(X_1,X_2) = g(X_1,X_1)g(X_2,X_2)-(g(X_1,X_2))^2 = a_{11} a_{22} - (a_{12})^2 = \det{\mathcal{A}}.
    \end{equation*}
    \item[$(5)$]
    \begin{equation*}
     K = \frac{g(R_{X_1 X_2}X_1,X_2)}{Q(X_1,X_2)} = \frac{g\left( \xi_1 X_1 + \xi_2 X_2,X_2 \right)}{\det{\mathcal{A}}^2} = \frac{a_{12}\xi_1 + a_{22}\xi_2}{\det{\mathcal{A}}^2} = \frac{a_{12}\xi_1 + a_{22}\xi_2}{\det{\mathcal{A}}^2}.
     \end{equation*}
     
\end{itemize}

\end{document}